\documentclass[leqno,twoside]{article} 
\usepackage[latin1]{inputenc}
\usepackage[T1]{fontenc}
\usepackage[english]{babel}
\usepackage{times}
\usepackage{amssymb} 
\usepackage{graphicx} 
\usepackage{tikz}
\usepackage{bm,dsfont}
\usepackage[margin=1.8cm,top=2.8cm]{geometry}
\usepackage{theorem}
\usepackage{multicol}
\usepackage{float}
\usepackage{hyperref}
\usepackage{cedya2015} 

\usepackage{amsmath}
\allowdisplaybreaks[4]

\newcommand*{\der}[2]{\frac{\partial #1}{\partial #2}}
\newcommand{\di}{\,\mathrm{d}}                              
\newcommand{\iin}{\;\text{in}\;} \newcommand{\oon}{\;\text{on}\;}
\newcommand{\deO}{{\partial\Omega}}
\newcommand*{\abs}[1]{\left|#1\right|}
\newcommand{\Tnorm}[1]{|||#1|||}
\newcommand*{\jmp}[1]{[\![#1]\!]}                     
\newcommand*{\mvl}[1]{\{\!\!\{#1\}\!\!\}}             

\newcommand{\bd}{{\mathbf d}}\newcommand{\be}{{\mathbf e}}\newcommand{\bn}{{\mathbf n}}
\newcommand{\bx}{{\mathbf x}}\newcommand{\bN}{{\mathbf N}}
\newcommand{\bT}{{\mathbf T}}\newcommand{\bV}{{\mathbf V}}\newcommand{\IR}{\mathbb{R}}
\newcommand{\bzeta}{{\boldsymbol\zeta}}      \newcommand{\bsigma}{{\boldsymbol \sigma}}
\newcommand{\btau}{{\boldsymbol\tau}}		     \newcommand{\bPhi}{{\boldsymbol\Phi}}      
\newcommand{\bzero}{{\mathbf 0}}
\newcommand{\calA}{{\mathcal A}}\newcommand{\calE}{{\mathcal E}}
\newcommand{\calF}{{\mathcal F}}\newcommand{\calT}{{\mathcal T}}
\newcommand{\Fh}{\calF_h}\newcommand{\Th}{{(\calT_h)}}
\newcommand{\deK}{{\partial K}}
\newcommand{\GD}{{\Gamma_D}}\newcommand{\GN}{{\Gamma_N}}\newcommand{\GR}{{\Gamma_R}}
\newcommand{\FT}{{\Fh^T}}\newcommand{\FO}{{\Fh^0}} 
\newcommand{\FD}{{\Fh^D}}\newcommand{\FN}{{\Fh^N}}\newcommand{\FR}{{\Fh^R}}
\newcommand{\hp}{_{hp}}
\newcommand{\vhp}{{v\hp}}
\newcommand\Vhp{v\hp}
\newcommand\Shp{\bsigma\hp}
\newcommand\hVhp{\widehat v\hp}
\newcommand\hShp{\widehat \bsigma\hp}
\newcommand{\Cstab}{C_{\mathrm{stab}}}
\newcommand*{\Norm}[1]{\left\|#1\right\|}
\newcommand{\DG}{_{\mathrm{DG}}}
\newcommand{\DGp}{_{\mathrm{DG^+}}}
\newcommand{\tht}{\vartheta} 
\newcommand{\Fspa}{{\Fh^{\text{space}}}}
\newcommand{\Ftime}{{\Fh^{\text{time}}}}
\newcommand{\mmbox}[1]{\fbox{\ensuremath{\displaystyle{ #1 }}}}	

\newcommand\heading{
\begin{flushright}
   \footnotesize 
   \textbf{\textsc{\textcolor{blue}{Proceedings of the XXIV Congress on Differential Equations and Applications}}}
   \\[0.3em]
\textbf{\textsc{\textcolor{blue}{XIV Congress on Applied Mathematics}}}
   \\[0.3em]
\textbf{\textrm{\textcolor{blue}{Cádiz, June 8-12, 2015, pp.~\thepage--\pageref{ultima-pagina:}}}}
   \\[1em]
\end{flushright}
} 

\begin{document}
\title{\heading \huge 
\textbf{Trefftz discontinuous Galerkin methods\\ on unstructured meshes for the wave equation}}


\author{A.~Moiola\thanks{Department of Mathematics and Statistics, University of Reading, Whiteknights PO Box 220, Reading RG6 6AX, UK. 
Email: a.moiola@reading.ac.uk}}

\markboth{~\hrulefill\ A.~Moiola}{Trefftz discontinuous Galerkin methods on unstructured meshes for the wave equation\hrulefill~}

\maketitle

\begin{Abstract}
We describe and analyse a space--time Trefftz discontinuous Galerkin method for the wave equation.
The method is defined for unstructured meshes whose internal faces need not be aligned to the space--time axes.
We show that the scheme is well-posed and dissipative, and we prove a priori error bounds for general Trefftz discrete spaces.
A concrete discretisation can be obtained using piecewise polynomials that satisfy the wave equation elementwise.


\vspace*{2ex}\par
\end{Abstract}


\begin{multicols}{2}



\section{Introduction}\label{s:Intro}

Finite element-type schemes whose discrete functions are piecewise solutions of the discretised PDE are often named \emph{Trefftz methods}.
When applied to time-dependent problems, Trefftz methods require the use of a space--time mesh, as opposed to the combination of a discretisation in space and a time-stepping.
In the context of linear wave propagation, Trefftz methods have been widely studied in time-harmonic regime (e.g.\ \cite{PVersion}), while only few recent works have been devoted to the time-domain acoustic and electromagnetic wave equations, chiefly \cite{EKSW2,SpaceTimeTDG,KSTW2014,PFT09,WTF14}.
Some possible advantages of Trefftz methods over standard ones are better orders of convergence; flexibility in the choice of basis functions; low dispersion; incorporation of wave propagation directions in the discrete space; adaptivity and local space--time mesh refinement.

Here we introduce a Trefftz discontinuous Galerkin (DG) method for the wave equation written as a first-order system. 
We consider unstructured meshes whose faces are not aligned to the space--time axes.
The metod proposed is an extension of those in \cite{SpaceTimeTDG,KSTW2014} to higher dimensions and more general meshes, and is closely related to those in \cite{EKSW2} (Maxwell's equations) and in \cite{MoRi05} (non-Trefftz, see Remark \ref{rem:MR} below).
After introducing the initial boundary value problem in \S\ref{s:IBVP}, some notation in \S\ref{s:Mesh} and the Trefftz DG formulation in \S\ref{s:TDG}, in \S\ref{s:WellP} we show that the Trefftz-DG formulation is well-posed for any Trefftz discrete space and we prove a priori error bounds in a DG norm.
In \S\ref{s:Duality}--\S\ref{s:Stability} we study the stability of an auxiliary problem in order to give a sufficient condition for an a priori error bound on the space--time $L^2$ norm of the Trefftz-DG error, and show that for a class of meshes this condition is verified.
These meshes can be obtained with the ``tent-pitching algorithm'' of \cite{EGSU05,FalkRichter1999,MoRi05}, and allow to treat the Trefftz-DG method as a semi-explicit scheme and to compute the solution element by element by solving a sequence of small local linear systems.
We also show that the proposed method is dissipative.
The analysis is carried out closely following that of \cite{PVersion,SpaceTimeTDG}; the duality argument used to control the $L^2$ norm of the error 
follows an idea of \cite{MOW99}.


\section{The initial boundary value problem}\label{s:IBVP}

We consider an initial boundary value problem (IBVP) posed on a space--time domain $Q=\Omega\times (0,T)$, where $\Omega\subset\R^n$ is a bounded Lipschitz polytope with outward unit normal $\bn_\Omega^x$, $n\in\N$ and $T>0$.
The boundary of $\Omega$ is divided in three parts, denoted $\GD$, $\GN$ and $\GR$, corresponding to Dirichlet, Neumann and Robin boundary conditions; one or two of them may be empty.
The wave equation IBVP reads as
\begin{align}\label{eq:IBVP}
\left\{\begin{aligned}
&\nabla v+\der{\bsigma}t = \bzero &&\iin Q,\\
&\nabla\cdot\bsigma+\frac{1}{c^2}\der{v}t = 0 &&\iin Q,
\\
&v(\cdot,0)=v_0, \quad \bsigma(\cdot,0)=\bsigma_0 &&\oon \Omega,
\\
&v=g_D, &&\oon \GD\times [0,T],
\\
&\bsigma\cdot\bn_\Omega^x=g_N, &&\oon \GN\times [0,T],
\\
&\frac\tht c v-\bsigma\cdot \bn_\Omega^x=g_R , &&\oon \GR\times [0,T].
\end{aligned}
\right.
\end{align}
Here $v_0,\bsigma_0,g_D,g_N,g_R$ are the problem data;
$c>0$ is the  wave speed, which is assumed to be piecewise constant and independent of $t$;
$\tht\in L^\infty(\GR\times[0,T])$ is an impedance parameter, which is assumed to be uniformly positive.
The gradient $\nabla$ and divergence $\nabla\cdot$ operators are meant in the space variable $\bx$ only.
%
The equations in \eqref{eq:IBVP} may be derived from the second-order scalar wave equation $-\Delta U+c^{-2}\der{{^2}}{t^2}U=0$ setting $v=\der Ut$ and $\bsigma=-\nabla U$.

\section{Space--time mesh and notation}\label{s:Mesh}

We impose a finite element mesh $\calT_h$ on the space-time domain $Q$.
We assume that all its elements are Lipschitz polytopes, so that each internal face $F=\partial K_1\cap \partial K_2$, for $K_1,K_2\in\calT_h$, with positive $n$-dimensional measure, is a subset of a hyperplane: 
$$F\subset\Pi_F:=\big\{(\bx,t)\in \IR^{n+1}:\; \bx\cdot\bn^x_F+t\, n^t_F=C_F\big\},$$ 
where $(\bn^x_F,n^t_F)$ is a unit vector in $\R^{n+1}$ and $C_F\in\R$. 
We assume that all the internal mesh faces belong to one of the following two classes:
\begin{align}
&\text{on an internal face } 
F=\partial K_1\cap \partial K_2 
\text{, either}
\nonumber\\
&\begin{cases}
c|\bn_F^x|< n_F^t & \text{and $F$ is called ``space-like'' face, or}\\
n_F^t=0 & \text{and $F$ is called ``time-like'' face.}
\end{cases}
\label{eq:HorVerFaces}
\end{align}
(See Remark \ref{rem:FullyUnstructured} for the extension to more general meshes.)
On space-like faces, by convention, we choose $n^t_F>0$, i.e.\ the normal unit vector $(\bn^x_F,n^t_F)$ points towards future time.
Intuitively, space-like faces are hypersurfaces lying under the characteristic cones and on which initial conditions might be imposed, while time-like faces are propagations in time of the faces of a mesh in space only.
We use the following notation:
\begin{align*}
\Fh&:=\bigcup\nolimits_{K\in\calT_h}\partial K \quad \text{(the mesh skeleton)},\\
\Fspa&:= \text{the union of the internal space-like faces,}\\
\Ftime&:= \text{the union of the internal time-like faces,}\\
\FO&:=\Omega\times\{0\},\hspace{12mm} \FT:=\Omega\times\{T\},\\ 
\FD&:=\GD\times[0,T],\qquad
\FN:=\GN\times[0,T],\\
\FR&:=\GR\times[0,T]. 
\end{align*}
(We will consider more specific meshes with $\Ftime=\emptyset$ in \S\ref{s:Stability} below.)
We denote the outward-pointing unit normal vector on $\deK$ by $(\bn^x_K,n^t_K)$.
For piecewise-continuous scalar ($w$) and vector ($\btau$) fields, we define averages $\mvl\cdot$, space normal jumps $\jmp{\cdot}_\bN$ and time (full) jumps $\jmp{\cdot}_t$ in the standard DG notation:
on $F=\deK_1\cap\deK_2$, $K_1,K_2\in\calT_h$,
\begin{align*}
\mvl{w}&:=\frac{w_{|_{K_1}}+w_{|_{K_2}}}2,\qquad
\mvl{\btau}:=\frac{\btau_{|_{K_1}}+\btau_{|_{K_2}}}2, \\
\jmp{w}_\bN&:= w_{|_{K_1}}\bn_{K_1}^x+w_{|_{K_2}} \bn_{K_2}^x,\\
\jmp{\btau}_\bN&:= \btau_{|_{K_1}}\cdot\bn_{K_1}^x+\btau_{|_{K_2}} \cdot\bn_{K_2}^x,\\
\jmp{w}_t&:= w_{|_{K_1}}n_{K_1}^t+w_{|_{K_2}} n_{K_2}^t=(w^--w^+)n^t_F,\\
\jmp{\btau}_t&:= \btau_{|_{K_1}} n_{K_1}^t+\btau_{|_{K_2}} n_{K_2}^t= (\btau^--\btau^+)n^t_F,
\end{align*}
where, on space-like faces, $w^-$ and $w^+$ denote the traces of the function $w$ from the adjacent elements at lower and higher times, respectively (and similarly for $\btau^\pm$).
We use the notation $\jmp{\cdot}_\bN$ to recall that $\jmp{w}_\bN$ is a vector field parallel to $\bn^x_F$ and $\jmp{\btau}_\bN$ is the jump of the normal component of $\btau$ only.
We recall also that in this formulas $|n^t_K|,|\bn^x_K|\le1$, and that one of the two might be zero; in particular $\jmp{w}_t=0$ and $\jmp{\btau}_t=\bzero$ on $\Ftime$.
The following identities can easily be shown:
%
\begin{align}\nonumber
&w^-\jmp{w}_t-\frac12 \jmp{w^2}_t=\frac1{2n^t_F}\jmp{w}_t^2,
\\
&\btau^-\jmp{\btau}_t-\frac12 \jmp{\btau^2}_t=\frac1{2n^t_F}|\jmp{\btau}_t|^2,
\label{eq:JumpIds}\\
&w^-\jmp{\btau}_\bN+\btau^-\cdot\jmp{w}_\bN-\jmp{w\btau}_\bN   
=\frac1{n^t_F}\jmp{w}_t\jmp{\btau}_\bN,
\nonumber\\
&\mvl{w}\jmp{\btau }_\bN+\mvl{\btau}\cdot\jmp{ w }_\bN=\jmp{w\btau}_\bN,
\nonumber
\end{align}
We assume that the mesh $\calT_h$ is chosen  so that the wave speed $c$ is constant in each element; as $c$ is independent of time, it may jump only across faces in $\Ftime$.

Finally, we define the local and the global \emph{Trefftz spaces}:
\begin{align*} 
\bT(K):=&\Big\{(w ,\btau )\in L^2(K)^{1+n}, \text{ s.t. }
\btau_{|_{\deK}}\in L^2(\deK)^n
,\\
&\der w t, \nabla\cdot \btau \in  L^2(K), \; \der \btau t, \nabla w \in L^2(K)^n,\\
&\underbrace{\nabla w+\der{\btau}t = \bzero, \;
\nabla\cdot\btau+\frac1{c^{2}}\der{w}t = 0}_{\text{(Trefftz property)}}
 \Big\}\;\forall K\in \calT_h,
\\
\bT(\calT_h):=&\Big\{(w ,\btau )\in L^2(Q)^{1+n},\text{ s.t. }\\
&\qquad(w_{|_K},\btau_{|_K})\in \bT(K) \;\forall K\in \calT_h\Big\}.
\end{align*}
The solution $(v,\bsigma)$ of IBVP \eqref{eq:IBVP} is assumed to belong to $\bT\Th$.



\section{The Trefftz-discontinuous Galerkin method}\label{s:TDG}

To obtain the DG formulation, we multiply the first two equations of \eqref{eq:IBVP} with test fields $\btau$ and $w$ and integrate by parts on a single mesh element $K\in\calT_h$:
\begin{align}
\label{eq:Elemental0}
-\int_K\bigg(v
\Big(
\nabla\cdot\btau+\frac1{c^2}\der wt 
\Big) +\bsigma\cdot\Big(
\der\btau t+\nabla w  
\Big)\bigg)\di V
\\
+\int_{\partial K} \bigg((v\,\btau + \bsigma\, w)\cdot\bn_K^x
+\Big(\bsigma\cdot\btau + \frac{v\,w}{c^2}\Big)\,n_K^t\bigg)\di S
=0.
\nonumber
\end{align}
We look for a discrete solution $(v\hp,\bsigma\hp)$ approximating $(v,\bsigma)$ in a finite-dimensional (arbitrary, at this stage) Trefftz space $\bV \Th\subset\bT\Th$.
We take the test field $(w,\btau)$ in the same space $\bV\Th$, thus the volume integral over $K$ in \eqref{eq:Elemental0} vanishes.
The traces of $v\hp$ and $\bsigma\hp$ on the mesh skeleton are approximated by the (single-valued) \emph{numerical fluxes} $\hVhp$ and $\hShp$, so that \eqref{eq:Elemental0} is rewritten as:
\begin{align}\label{eq:Elemental1}
\int_{\partial K} 
\hVhp\Big(\btau \cdot\bn_K^x+\frac{w}{c^2} n_K^t\Big)
+\hShp\cdot\big(w\bn_K^x+\btau\,n_K^t\big)
\di S=0.
\end{align}
We choose to define the numerical fluxes as:
\begin{align*}
&\hVhp := \hspace{27mm} \hShp:=
\\
&\begin{cases}
\Vhp^- \\
\Vhp \\
v_0\\ 
\mvl{\Vhp}+\beta \jmp{\Shp}_\bN \\
g_D\\
v\hp\!+\!\beta(\bsigma\hp\cdot\bn_\Omega^x\!-\!g_N)\\
(1-\delta)v\hp\\
\quad+\frac{\delta c}\tht (\bsigma\hp\cdot\bn_\Omega^x\!+\!g_R)
\end{cases}\hspace{-4mm}
\begin{cases}
\Shp^- & \oon \Fspa\!,\\
\Shp &\oon \FT,\\
\bsigma_0 &\oon \FO,\\
\mvl{\Shp}+\alpha \jmp{\Vhp}_\bN \hspace{-2mm}&\oon \Ftime\!,\\
\Shp\!+\!
\alpha(v\hp\!-\!g_D)\bn_\Omega^x\hspace{-2mm}&\oon \FD,\\
g_N\bn_\Omega^x&\oon \FN,\\
(1\!-\!\delta)(\tht\frac {v\hp}c\!-\!g_R)\bn_\Omega^x\hspace{-2mm}&\oon \FR.\\
\qquad +\delta\bsigma\hp
\end{cases}
\end{align*}
The mesh-dependent parameters $\alpha\in L^\infty(\Ftime\cup\FD)$,  $\beta\in L^\infty(\Ftime\cup\FN)$ and  $\delta\in L^\infty(\FR)$ may be used to tune the method, e.g.\ to deal with locally-refined meshes.

\begin{figure*}[hbt]
\noindent\fbox{\parbox{0.98\textwidth}
{
\begin{align}\label{eq:TDG}
\text{Seek\;} (\Vhp,\Shp&)\in \bV (\calT_h) \text{ such that }\quad
\calA(\Vhp,\Shp; w ,\btau )=\ell( w ,\btau )\quad 
\forall ( w ,\btau )\in \bV (\calT_h), \text{ where}
\\
\calA(\Vhp,\Shp; w ,\btau ):=&
\int_{\Fspa}\big(c^{-2}\Vhp^-\jmp{w}_t+\Shp^-\cdot\jmp{\btau}_t+\Vhp^-\jmp{\btau}_\bN+\Shp^-\cdot\jmp{w}_\bN\big)\di S
+\int_\FT (c^{-2}\Vhp  w +\Shp \cdot\btau )\di \bx
\nonumber\\
&+\int_\Ftime \big( \mvl{\Vhp}\jmp{\btau }_\bN+\mvl{\Shp}\cdot\jmp{ w }_\bN
+\alpha\jmp{\Vhp}_\bN\cdot\jmp{ w }_\bN+ \beta\jmp{\Shp}_\bN\jmp{\btau }_\bN
\big)\di S
\nonumber\\
&+\int_\FD \big(\bsigma\cdot\bn_\Omega^x\, w +\alpha \Vhp w   \big) \di S+\int_\FN\big(v\hp(\btau\cdot\bn_\Omega^x)
+\beta(\bsigma\cdot\bn_\Omega^x)(\btau\cdot\bn_\Omega^x)\big)\di S
\nonumber\\
&+\int_\FR\Big(
\frac{(1-\delta)\tht}c v\hp w+(1-\delta)v\hp(\btau\cdot\bn_\Omega^x)
+\delta(\bsigma\hp\cdot\bn_\Omega^x) w+\frac{\delta c}\tht(\bsigma\hp\cdot\bn_\Omega^x) (\btau\cdot\bn_\Omega^x)
\Big)\di S,
\nonumber\\
\ell( w ,\btau ):=&
\int_\FO ( c^{-2}v_0 w  +\bsigma_0\cdot \btau )\di \bx
\nonumber\\&
+\int_\FD g_D\big(\alpha  w -\btau\cdot\bn_\Omega^x\big)\di S
+\int_\FN g_N \big(\beta\,\btau\cdot\bn_\Omega^x-w\big)\di S
+\int_\FR g_R \Big((1-\delta)w-\frac{\delta c}\tht\,\btau\cdot\bn_\Omega^x\Big)\di S.
\nonumber
\end{align}\vspace{-2mm}
}}
\end{figure*}

The fluxes are \emph{consistent}, in the sense that they coincide with the exact solution $(v,\bsigma)$ of the IBVP \eqref{eq:IBVP} if they are applied to $(v,\bsigma)$ itself, which satisfies the boundary conditions and has no jumps across mesh faces.
Moreover, the fluxes satisfy the boundary conditions, e.g.\ $\frac\tht c \hVhp-\hShp\cdot\bn_\Omega^x=g_R$ on $\FR$.
The numerical fluxes can be understood as upwind fluxes on the space-like faces and centred fluxes with jump penalisation on the time-like ones.

Summing the elemental TDG equation \eqref{eq:Elemental1} over the elements $K\in \calT_h$, with the fluxes defined above, we obtain the Trefftz-DG variational formulation displayed in box \eqref{eq:TDG}.


Method \eqref{eq:TDG} appears as a  formulation over the whole space--time domain $Q$, which would lead to an unconvenient global linear system coupling all elements.
However, if the mesh is suitably designed, the solution can be computed by solving a sequence of smaller local problem.
A first possibility is to partition the time interval $[0,T]$ into subintervals and solve sequentially for the corresponding time-slabs $\Omega\times(t_{j-1},t_j)$, see \cite{EKSW2,KSTW2014}; this corresponds to an implicit method and allows local mesh refinement.
A slightly more complicated, but potentially much more efficient version is to solve for small patches of elements, localised in space and time, in the spirit of the semi-explicit ``tent-pitching'' algorithm of \cite{EGSU05,FalkRichter1999,MoRi05}. 
If no time-like faces are present in the mesh, the solution can be found by solving a small system for each element, see \S\ref{s:Stability} below.
If the same mesh is used, all these variants are 
equivalent, in the sense that the discrete solutions $(\vhp,\Shp)$ coincide.

\begin{remark}\label{rem:FullyUnstructured}
One could easily extend the formulation weakening assumption \eqref{eq:HorVerFaces} to allow more general time-like faces with $c|\bn^x_F|>n^t_F$, i.e.\ not aligned to the time-axis.
Choosing the numerical fluxes as above, one obtains a formulation similar to \eqref{eq:TDG} with additional terms on $\Ftime$ containing $\jmp{w}_t$ and $\jmp{\btau}_t$, which do not vanish in this setting.
It is then easy to prove the coercivity of the new bilinear form in a slight modification of the DG norm introduced below.
However, the bilinear form will contain the term $\int_\Ftime\mvl{\Shp}\cdot\jmp{\btau}_t\di S$, featuring the full jump of $\btau$ (as opposed to the normal jump only), which, in dimension $n>1$, is not easily controlled by the same DG norms.
\end{remark}

\begin{remark}\label{rem:MR}
Formulation \eqref{eq:TDG} can be seen in the framework of DG methods for general first-order hyperbolic systems developed in \cite{MoRi05}, which considers standard discontinuous piecewise-polynomial spaces.
The choice of the numerical fluxes in the interior elements correspond to the choice of a suitable decomposition of the block matrix 
$\mathsf M=(\begin{smallmatrix}n^t_K c^{-2} &(\bn^x_K)^\top\\ \bn^t_K & n^t_K \mathrm{Id}_n \end{smallmatrix})$ 
defined on $\deK$ for all $K\in \calT_h$.
Here $\mathrm{Id}_n$ is the identity matrix in $\IR^{n\times n}$, ${}^\top$ denotes vector transposition, and $\bn^x_K$ is thought as a column vector.
The choice we have done in this section corresponds to the decomposition:
\setlength{\tabcolsep}{2pt}
\begin{tabular}{l|l|l}
on & $\mathsf M^+ =$ & $\mathsf M^-=$ \\
\hline
$\!\partial K\cap\Fspa\cap\{n^t_K\!>\!0\}$&
$\mathsf M$&
$\bzero$
\\
$\!\partial K\cap\Fspa\cap\{n^t_K\!<\!0\}$&$\bzero$&
$\mathsf M$
\\
$\!\partial K\cap\Ftime$ 
&$(\begin{smallmatrix}\alpha &\frac12(\bn^x_K)^\top\\ 
\frac12\bn^x_K & \beta\bn^x_K\otimes \bn^x_K\end{smallmatrix})$&
$(\begin{smallmatrix}-\alpha &\frac12(\bn^x_K)^\top\\ 
\frac12\bn^x_K & -\beta\bn^x_K\otimes \bn^x_K\end{smallmatrix})$
\end{tabular}

The conditions $\mathsf M^+ +\mathsf M^-=\mathsf M $,
$\ker(\mathsf M^+ -\mathsf M^-)=\ker(\mathsf M)$
and $\mathsf M^+_{|_{\deK_1}} +\mathsf M^-_{|_{\deK_2}}=\bzero$ on $\deK_1\cap\deK_2$
are automatically satisfied.
Moreover, $\mathsf M^+\ge0$ and $\mathsf M^-\le 0$ hold true if and only if $\alpha\beta\ge1/4$.
The boundary terms in \eqref{eq:TDG} and in \cite[\S6]{MoRi05} coincide (apart from a different sign convention) if our flux parameters and their boundary coefficients $Q$ and $\sigma$ are chosen so that 
$\alpha=\sigma$ on $\GD$, $\beta=1/\sigma$ on $\GN$,  $\delta=(1+Q)/2$ and $\tht/c=\sigma(1+Q)/(1-Q)$ on $\GR$.
\end{remark}

\section{A priori error analysis}\label{s:Analysis}

\subsection{Definitions and assumptions}\label{s:Assumptions}

We prove the well-posedness and the stability of the Trefftz-DG formulation \eqref{eq:TDG} under the assumption that the flux parameters $\alpha$, $\beta$ and $\delta$ are uniformly positive in their domains of definition and that $\Norm{\delta}_{L^\infty(\FR)}<1$.
We introduce a piecewise-constant function $\gamma$ defined on $\Fspa$, measuring how close to characteristic cones the space-like mesh faces are:
\begin{align}\label{eq:gamma}
\gamma:= \frac{c|\bn^x_F|}{n^t_F}\;\oon F
\subset\Fspa, \qquad \gamma:=0\;\oon \FO\cup\FT,
\end{align}
from which, recalling assumption \eqref{eq:HorVerFaces}, $\gamma\in(0,1)$ and 
\begin{align}\label{eq:JumpIneq}
&
\abs{\jmp{w}_\bN}\le \frac \gamma c\abs{\jmp{w}_t},\quad
\abs{\jmp{\btau}_\bN}\le \frac \gamma c\abs{\jmp{\btau}_t} \quad\oon \Fspa.
\end{align}
%
We define two mesh- and flux-dependent norms on $\bT\Th$:
\begin{align*}
&\Tnorm{(w ,\btau )}^2\DG:=
\frac12 \Norm{\Big(\frac{1-\gamma}{n^t_F}\Big)^{1/2}c^{-1}\jmp{w }_t}_{L^2(\Fspa)}^2\\
&\quad
+\frac12\Norm{\Big(\frac{1-\gamma}{n^t_F}\Big)^{1/2}  
\jmp{\btau }_t\rule{0pt}{3mm}}_{L^2(\Fspa)^n}^2 \\
&\quad
+\frac12\Norm{c^{-1}w }^2_{L^2(\FO\cup\FT)} 
+\frac12\Norm{\btau \rule{0pt}{3mm}}^2_{L^2(\FO\cup\FT)^n} \\
&\quad
+\Norm{\alpha^{1/2}\jmp{w }_\bN}_{L^2(\Ftime)^n}^2
+\Norm{\beta^{1/2}\jmp{\btau }_\bN}_{L^2(\Ftime)}^2\\
&\quad+\Norm{\alpha^{1/2}w }_{L^2(\FD)}^2
+\Norm{\beta^{1/2} \btau\cdot\bn_\Omega^x }_{L^2(\FN)}^2\\
&\quad+\Norm{\Big(\frac{(1-\delta)\tht}c\Big)^{1/2} w }_{L^2(\FR)}^2
+\Norm{\Big(\frac{\delta c}\tht\Big)^{1/2} \btau\cdot\bn_\Omega^x }_{L^2(\FR)}^2\!;
\\
&\Tnorm{(w ,\btau )}^2\DGp:=
\Tnorm{(w ,\btau )}^2\DG\\
&\quad
+\Norm{\Big(\frac{n^t_F}{1-\gamma}\Big)^{1/2} 
c^{-1}w^-}_{L^2(\Fspa)}^2\\
&\quad
+\Norm{\Big(\frac{n^t_F}{1-\gamma}\Big)^{1/2}
\btau^-}_{L^2(\Fspa)^n}^2  \\
&\quad
+\Norm{
\beta^{-1/2}\mvl{w }}_{L^2(\Ftime)}^2
+\Norm{
\alpha^{-1/2}\mvl{\btau }}_{L^2(\Ftime)^n}^2\\
&\quad
+\Norm{\alpha^{-1/2}\btau\cdot\bn_\Omega^x }_{L^2(\FD)}^2
+\Norm{\beta^{-1/2} w }_{L^2(\FN)}^2.
\end{align*}
(The factors $(1-\gamma)^{\pm1/2}$ may be dropped from all terms in the two norms without modifying the analysis in the following sections.)
We note that these are only \textit{seminorms} on broken Sobolev spaces defined on the mesh $\calT_h$, 
but are norms on $\bT\Th$: indeed $\Tnorm{(w,\btau)}=0$ for $(w,\btau)\in\bT\Th$ implies that  $(w,\btau)$ is solution of the IBVP \eqref{eq:IBVP} with zero initial and boundary conditions, so $(w,\btau)=(0,\bzero)$ by the well-posedness of the IBVP itself (see \cite[Lemma~4.1]{SpaceTimeTDG}).


\subsection{Well-posedness}\label{s:WellP}

We first note that for all Trefftz field $(w,\btau)\in\bT\Th$
$$\int_{\deK}
\bigg(w\btau\cdot\bn^x_K+\frac12\Big(\frac{w^2}{c^2}+|\btau|^2\Big)n^t_K\bigg)\di S=0\quad \forall K\in\calT_h,$$
which follows from integration by parts and the Trefftz property.
Subtracting these terms from the bilinear form $\calA$, 
using the jump identities \eqref{eq:JumpIds}, the inequalities \eqref{eq:JumpIneq}, the definition of $\gamma$ in \eqref{eq:gamma}, and the weighted Cauchy--Schwarz inequality, we show that the form $\calA$ is coercive in $\Tnorm{\cdot}\DG$ norm with unit constant:
\begin{align*}
&\calA(w,\btau; w ,\btau )\\
&=
\calA(w,\btau; w ,\btau )
\\&\qquad
-\sum_{K\in\calT_h}\int_{\deK}
\bigg(w\btau\cdot\bn^x_K+\frac{1}2\Big(\frac{w^2}{c^2}+|\btau|^2\Big)n^t_K\bigg)\di S
\\
&\overset{\eqref{eq:TDG}}=\int_{\Fspa}\big(c^{-2}w^-\jmp{w}_t+\btau^-\cdot\jmp{\btau}_t
+w^-\jmp{\btau}_\bN+\btau^-\cdot\jmp{w}_\bN
\\&\qquad
-\jmp{w\btau}_\bN-\frac12\jmp{c^{-2}w^2+|\btau|^2}_t\big)\di \bx
\\&
+\frac12\int_\FT (c^{-2}w^2  +\abs{\btau}^2 )\di \bx
+\frac12\int_\FO(c^{-2}w^2  +\abs{\btau}^2 )\di \bx
\\&
+\int_\Ftime \Big( \mvl{w}\jmp{\btau }_\bN+\mvl{\btau}\cdot\jmp{ w }_\bN
\\&\qquad
+\alpha|\jmp{w}_\bN|^2+ \beta\jmp{\btau}_\bN^2
-\jmp{w\btau}_\bN
\Big)\di S
\\&
+\int_\FD \alpha  w^2    \di S
+\int_\FN \beta(\btau\cdot\bn_\Omega^x)(\btau\cdot\bn_\Omega^x)\di S
\\&
+\int_\FR\Big(\frac{(1-\delta)\tht}c w^2
+\frac{\delta c}\tht(\btau\cdot\bn_\Omega^x)^2
\Big)\di S
\\&
\overset{\eqref{eq:JumpIds}}=
\int_{\Fspa}\bigg(
\frac1{2n^t_F}(c^{-2}\jmp{w}_t^2+|\jmp{\btau}_t|^2)+\frac1{n^t_F}\jmp{w}_t\jmp{\btau}_\bN
\bigg)\di \bx
\\&
+\frac12\Norm{c^{-1}w }^2_{L^2(\FO\cup\FT)} 
+\frac12\Norm{\btau \rule{0pt}{3mm}}^2_{L^2(\FO\cup\FT)^n}
\\&
+\Norm{\alpha^{1/2}\jmp{w}_\bN}^2_{L^2(\Ftime)}
+\Norm{\beta^{1/2}\jmp{\btau}_\bN}^2_{L^2(\Ftime)}
\\&
+\Norm{\alpha^{1/2}w }_{L^2(\FD)}^2
+\Norm{\beta^{1/2} \btau\cdot\bn_\Omega^x }_{L^2(\FN)}^2
\\&
+\Norm{\Big(\frac{(1-\delta)\tht}c\Big)^{1/2} w }_{L^2(\FR)}^2
+\Norm{\Big(\frac{\delta c}\tht\Big)^{1/2} \btau\cdot\bn_\Omega^x }_{L^2(\FR)}^2
\\&
\overset{\eqref{eq:gamma},\eqref{eq:JumpIneq}}
\ge
\Tnorm{(w,\btau)}\DG^2
\qquad\qquad \forall (w,\btau)\in \bT\Th.
\end{align*}
Using again the Cauchy--Schwarz inequality, the bounds on the jumps \eqref{eq:JumpIneq} and $\gamma<1$, we have the following continuity estimate for the bilinear form $\calA$: for all $(v,\bsigma),(w ,\btau )\in\bT\Th$
\begin{align}
\nonumber
&\abs{\calA(v,\bsigma; w ,\btau )}\le C_c\Tnorm{(v,\bsigma)}\DGp\Tnorm{(w,\btau)}\DG,
\quad \text{where}
\\&
C_c:=
\begin{cases}
2, &\FR=\emptyset,\\
2\max\Big\{\Norm{\frac{1-\delta}\delta}_{L^\infty(\FR)}^{1/2},
\big\|\frac\delta{1-\delta}\big\|_{L^\infty(\FR)}^{1/2}\Big\}
&\FR\ne\emptyset.
\end{cases}
\label{eq:Continuity}
\end{align}
Note that $C_c\ge2$ and that the minimal value $C_c=2$ is obtained for $\delta=1/2$.
Also the linear functional $\ell$ is continuous:
\begin{align*}
\abs{\ell(w,\btau)}
\le&\Big(2\Norm{c^{-1}v_0}_{L^2(\FO)}^2+2\Norm{\bsigma_0}_{L^2(\FO)}^2
\\&
+2\Norm{\alpha^{1/2}g_D}_{L^2(\FD)}^2+2\Norm{\beta^{1/2}g_N}_{L^2(\FN)}^2
\\&
+
\Norm{(c/\tht)^{1/2}g_R}_{L^2(\FR)}^2
\Big)^{1/2} \Tnorm{(w,\btau)}\DGp;
\end{align*}
if $g_D=g_N=0$ (or the corresponding parts $\FD,\FN$ of the boundary are empty) then 
the $\Tnorm{(w,\btau)}\DGp$ norm at the right-hand side can be substituted by $\Tnorm{(w,\btau)}\DG$.

Combining coercivity and continuity, C\'ea lemma gives well-posedness and quasi-optimality of the Trefftz-DG formulation irrespectively of the discrete Trefftz space $\bV\Th$ chosen.

\begin{theorem}\label{th:QO}
The variational problem \eqref{eq:TDG} admits a unique solution $(v\hp,\bsigma\hp)\in \bV\Th$.
It satisfies:
\begin{align*}
\mmbox{\begin{aligned}&\Tnorm{(v-\Vhp,\bsigma-\Shp)}\DG\\
&\qquad\le  (1+C_c) \inf_{(w,\btau)\in\bV\Th} \Tnorm{(v-w,\bsigma-\btau)}\DGp,
\end{aligned}}
\end{align*}
with $C_c$ as in \eqref{eq:Continuity}.
Moreover, if $g_D=g_N=0$ then
\begin{align*}
&\Tnorm{(v\hp,\bsigma\hp)}\DG\le
\\& 
\bigg(2\Norm{\frac1c 
v_0}_{L^2(\FO)}^2+2\Norm{\bsigma_0}_{L^2(\FO)}^2
+\Norm{\frac{c^{1/2}}{\tht^{1/2}}g_R}_{L^2(\FR)}^2
\bigg)^{1/2}\!.
\end{align*}
\end{theorem}
In the next sections we control the $L^2(Q)$ norm of the error with a duality-type argument inspired by \cite[Theorem~3.1]{MOW99}, based on an auxiliary problem and on energy estimates.


\subsection{Auxiliary problem and error bounds in \texorpdfstring{$L^2(Q)$}{L2(Q)}}\label{s:Duality}
Consider the auxiliary inhomogeneous IBVP
\begin{align}\label{eq:zzIBVP}
\left\{\begin{aligned}
&\nabla z+ \partial\bzeta/\partial t 
= \bPhi &&\iin Q,\\
&\nabla\cdot\bzeta+c^{-2} \,\partial z/\partial t 
= \psi &&\iin Q,\\
&z(\cdot,0)=0, \quad \bzeta(\cdot,0)=\bzero\qquad &&\oon \Omega,\\
&z=0 &&\oon \GD\times I,\\
&\bzeta\cdot\bn_\Omega^x=0 &&\oon \GN\times I,\\
&\frac\tht c z-\bzeta\cdot\bn_\Omega^x=0 &&\oon \GR\times I,
\end{aligned}\right.\end{align}
with data $\psi\in L^2(Q),\bPhi\in L^2(Q)^n$.
In the next proposition we assume that there exists a positive constant $\Cstab$, depending on the domain $Q$, the mesh $\calT_h$ (thus on $\gamma$) and the parameters $c,\tht,\alpha,\beta,\delta$, such that for all 
$\psi\in L^2(Q),\bPhi\in L^2(Q)^n$ the solution $(z,\bzeta)$ of \eqref{eq:zzIBVP} satisfies the stability bound:
\begin{align}\label{eq:DualStability}
&
2\Norm{\Big(\frac{(1+\gamma^2)n^t_F}{1-\gamma}\Big)^{1/2}\frac zc}^2_{L^2(\Fspa\cup\FT)}
\\
&+2\Norm{\Big(\frac{(1+\gamma^2)n^t_F}{1-\gamma}\Big)^{1/2}\bzeta}^2_{L^2(\Fspa\cup\FT)^n}
\nonumber\\&
+\Norm{\alpha^{-1/2}\bzeta\cdot\bn^x_F}^2_{L^2(\Ftime\cup\FD)}
+\Norm{\beta^{-1/2}z}^2_{L^2(\Ftime\cup\FN)}
\nonumber\\&
+\Norm{\Big(\frac c{(1-\delta)\tht}\Big)^{1/2}\bzeta\cdot\bn^x_\Omega}^2_{L^2(\FR)}
+\Norm{\Big(\frac \tht{\delta c}\Big)^{1/2}z}^2_{L^2(\FR)}
\nonumber\\
&\le
\Cstab^2\Big(\Norm{\bPhi}^2_{L^2(Q)^n}+\Norm{c\psi}^2_{L^2(Q)}\Big).
\nonumber
\end{align}
Conditions under which bound \eqref{eq:DualStability} holds are given in \S\ref{s:Stability}.

\begin{proposition}\label{prop:Duality}
If bound \eqref{eq:DualStability} is satisfied, then for all Trefftz fields $(w,\btau)\in\bT\Th$
\begin{align}
{\Big(\Norm{c^{-1}w}_{L^2(Q)}^2+\Norm{\btau}_{L^2(Q)^n}^2\Big)^{1/2}
\le \Cstab \Tnorm{(w,\btau)}\DG.}
\label{eq:zzDualityBound}
\end{align}
\end{proposition}

\textit{Proof.}\;
We first prove the vanishing of certain jumps across mesh faces for the solution $(z,\bzeta)$ of the inhomogeneous auxiliary problem \eqref{eq:zzIBVP}:
$\jmp{z}_t$ and $\jmp{\bzeta}_t$ on $\Fspa$ and of $\jmp{z}_\bN$ and $\jmp{\bzeta}_\bN$ on $\Ftime$.
%
%
%
Given a hyperplane $\Pi=\{\bn^x_\Pi\cdot \bx+n^t_\Pi\,t=C_\Pi\}$, 
denote the scalar jump of a piecewise continuous function $f$ as 
$\jmp{f}_\Pi:=f_{|\{\bn^x_\Pi\cdot \bx+n^t_\Pi\,t>C_\Pi\}}-f_{|\{\bn^x_\Pi\cdot \bx+n^t_\Pi\,t<C_\Pi\}}$.
From \eqref{eq:zzIBVP}, the fields $(\bzeta, c^{-2}z)$ and $(z\be_j,\bzeta_j)$, $1\le j\le n$ ($\be_j$ denoting the standard basis elements of $\R^n$), are in $H(\mathrm{div}_{x,t},Q)$, thus their \textit{normal} jumps vanish across any space--time Lipschitz interface in $Q$ and in particular across $\Pi$:
$$\jmp{\bzeta\cdot\bn^x_\Pi + c^{-2}z n^t_\Pi}_\Pi=\jmp{z(\bn^x_\Pi)_j + \bzeta_j n^t_\Pi}_\Pi=0
\qquad 1\le j\le n.$$
Thus, on time-like faces $n^t_\Pi=0$, the jump of $z$ and the normal jump of $\bzeta$ vanish.
On constant-time faces ($\bn^x_\Pi=\bzero, n^t_\Pi=\pm1$) all jumps vanish.
On other planes, $n^t_\Pi\ne0$ and $|\bn^x_\Pi|\ne 0$, thus
$$(-c^2/n^t_\Pi)\jmp{\bzeta\cdot\bn^x_\Pi }_\Pi= \jmp{z}_\Pi
=(-n^t_\Pi/|\bn^x_\Pi|^2)\jmp{\bzeta\cdot\bn^x_\Pi }_\Pi.$$
If $n^t_\Pi/|\bn^x_\Pi|\ne c$ then we have immediately $\jmp{z}_\Pi=\jmp{\bzeta\cdot\bn^x_\Pi }_\Pi=0$  and from above, $\jmp{\bzeta_j}_\Pi=0$ for all $1\le j\le n$. 
Assumption \eqref{eq:HorVerFaces} guarantees that $n^t_\Pi/|\bn^x_\Pi|> c$ on $\Fspa$, so we conclude that all jumps vanish.
(We have simply shown that the discontinuities of solutions of the first-order wave equation with source term in $L^2(Q)^{n+1}$ propagate along characteristics.)

Since the $L^2(Q)$ norm of  $(w,\btau)\in\bT\Th$ can be computed as
\begin{align}
\label{eq:L2sup}
&\big(\Norm{c^{-1}w}_{L^2(Q)}^2+\Norm{\btau}_{L^2(Q)}^2\big)^{1/2}
\\&=
\sup_{(\bPhi,\psi)\in L^2(Q)^{n+1}}\frac{\iint_Q(w\psi+\btau\cdot\bPhi)\di x\di t}{\big(\Norm{\bPhi}^2_{L^2(Q)^n}+\Norm{c\psi}^2_{L^2(Q)}\big)^{1/2}},
\nonumber
\end{align}
we now take the scalar product of $(w,\btau)$ with the source terms $(\psi,\bPhi)$ of problem \eqref{eq:zzIBVP} and integrate by parts in each element:
\begin{align*}
&\int_Q(w\psi+\btau\cdot\bzeta)\di V\\
\overset{\eqref{eq:zzIBVP}}=&
\sum_{K\in\calT_h}\int_K \bigg(w\nabla\cdot\bzeta+\frac{w}{c^2}\der zt
+\btau\cdot\nabla z+\btau\cdot \der\bzeta t\bigg)\di V\\
=&\sum_{K\in\calT_h}\int_\deK \Big(
w\bzeta\cdot\bn^x_K+\btau\cdot\bn^x_K z+\frac{wz n^t_K}{c^2}+\btau\cdot \bzeta n^t_K\Big)\di S
\\
=&\int_\Fspa\underbrace{\jmp{w\bzeta+\btau z}_\bN+\jmp{{c^{-2}}wz +\btau\cdot \bzeta}_t}_{
\le c^{-1}|\jmp{w}_t|(\gamma |\bzeta|+c^{-1}|z|)+|\jmp{\btau}_t|(\gamma c^{-1}|z|+|\bzeta|)}
\di S
\\&
+\int_{\FT}\Big(\frac{wz}{c^2} +\btau\cdot \bzeta\Big)\di S
-\int_{\FO}\Big(\frac w{c^2}\underbrace{z}_{=0} 
+\btau\cdot \underbrace{\bzeta}_{=\bzero}\Big)\di S
\\&
+\int_{\Ftime}\underbrace{\jmp{w\bzeta+\btau z}_\bN}_{=\jmp{w}_\bN\cdot\bzeta+\jmp{\btau}_\bN z}\di S
\\&
+\int_{\FD\cup\FN\cup\FR}\big(w\underbrace{\bzeta\cdot\bn^x_\Omega}_{=0 \oon\FN}
+\btau\cdot\bn^x_\Omega \underbrace{z}_{=0\oon \FD}\big)\di S
\\
\le&
\Tnorm{(w,\btau)}\DG \cdot\bigg(
2\Norm{\Big(\frac{(1+\gamma^2)n^t_F}{1-\gamma}\Big)^{1/2}\frac zc}^2_{L^2(\Fspa)}
\\&
+2\Norm{\Big(\frac{(1+\gamma^2)n^t_F}{1-\gamma}\Big)^{1/2}\bzeta}^2_{L^2(\Fspa)^n}
\\&
+2\Norm{c^{-1}z}^2_{L^2(\FT)}
+2\Norm{\bzeta}^2_{L^2(\FT)^n}
\\&
+\Norm{\alpha^{-1/2}\bzeta\cdot\bn^x_F}^2_{L^2(\Ftime\cup\FD)}
+\Norm{\beta^{-1/2}z}^2_{L^2(\Ftime\cup\FN)}
\\&
+\Norm{\Big(\frac c{(1-\delta)\tht}\Big)^{1/2}\bzeta\cdot\bn^x_\Omega}^2_{L^2(\FR)}
\\&
+\Norm{\Big(\frac \tht{\delta c}\Big)^{1/2}z}^2_{L^2(\FR)}
\bigg)^{1/2}
\\
\overset{\eqref{eq:DualStability}}\le& \Cstab\Tnorm{(w,\btau)}\DG
\big(\Norm{\bPhi}^2_{L^2(Q)^n}+\Norm{c\psi}^2_{L^2(Q)}\big)^{1/2}.
\end{align*}
Inserting this bound in the expansion \eqref{eq:L2sup} of the $L^2(Q)$ norm of $(w,\btau)$, we obtain the assertion \eqref{eq:zzDualityBound}.
\hfill $\square$

From Proposition \ref{prop:Duality} and Theorem \ref{th:QO}, it follows that the $L^2(Q)$ norm of the Trefftz-DG error is controlled by the $\Tnorm{\cdot}\DGp$ norm of the best-approximation error: 
if bound \eqref{eq:DualStability} is verified,
\begin{align}\label{eq:L2Error}
\mmbox{\begin{aligned}
&\Big(\Norm{c^{-1}(v-\Vhp)}_{L^2(Q)}^2+\Norm{\bsigma-\Shp}_{L^2(Q)^n}^2\Big)^{1/2}\\
&\le  \Cstab(1\!+\!C_c)\!\! \inf_{(w,\btau)\in\bV\Th} \Tnorm{(v-w,\bsigma-\btau)}\DGp,
\end{aligned}}
\end{align}

In \cite{SpaceTimeTDG}, bound \eqref{eq:DualStability} was proven in one space dimension on meshes made of rectangular elements aligned to the space-time axes and $\Cstab$ was computed explicitly.
Two proofs were given.
One of them (Appendix~A of \cite{SpaceTimeTDG}) relies on the use of the exact value of $(z,\bzeta)$ in $Q$ computed with Duhamel's principle, 
and can not be easily extended to general domains in higher dimensions as it require a suitable periodic extension of the IBVP \eqref{eq:zzIBVP} to $\R^n\times(0,T)$.
The second proof (Lemma~4.9 of \cite{SpaceTimeTDG}) uses an energy argument to control the traces on space-like faces in \eqref{eq:DualStability} and an integration by parts trick to bound the traces on time-like faces.
In higher dimensions, the energy argument carries over, while the traces on time-like faces are harder to control.
In \S\ref{s:Stability} we prove the stability estimate \eqref{eq:DualStability} in any dimension, under two additional assumptions to get around the need to control traces on time-like faces.
We first discuss in the next section some energy identities which will be useful later.

\subsection{Energy identities and estimates}\label{s:Energy}

We call ``space-like interface'' a graph hypersurface 
$$\Sigma=\big\{(\bx,f_\Sigma(\bx)):\; \bx\in \Omega\big\}\subset\overline Q$$ 
where $f_\Sigma:\overline \Omega\to[0,T]$ is a Lipschitz-continuous function whose Lipschitz constant is smaller than $1/c$.
Each space-like mesh face in $\Fspa$ is subset of a space-like interface $\Sigma$.
The future-pointing unit normal vector is defined almost everywhere on $\Sigma$ and denoted by $(\bn^x_\Sigma,n^t_\Sigma)$.

For sufficiently smooth scalar and vector fields $(w,\btau)$, define their \emph{energy} on a space-like interface as
$$\calE(\Sigma;w,\btau):=\int_\Sigma \bigg(
w\,\btau\cdot \bn^x_\Sigma+\frac12\Big(\frac{w^2}{c^2}+|\btau|^2\Big)n^t_\Sigma
\bigg)\di S.
$$
The energy on constant-time, or ``flat'', space-like interfaces are denoted by $\calE(t;w,\btau):=\calE(\overline\Omega\times\{t\};w,\btau)$, for $0\le t\le T$.

For two space-like interfaces $\Sigma_1,\Sigma_2$ with $f_{\Sigma_1}\le f_{\Sigma_2}$ in $\overline \Omega$, we denote the volume between them and its lateral boundary as
\begin{align*}
Q_{\Sigma_1,\Sigma_2}&:=\{(\bx,t): \;\bx\in\Omega,\; f_{\Sigma_1}(\bx)<t<f_{\Sigma_2}(\bx)\},\\
\Gamma_{\Sigma_1,\Sigma_2}&:=\{(\bx,t):\; \bx\in\deO,\; f_{\Sigma_1}(\bx)\le t\le f_{\Sigma_2}(\bx)\}.
\end{align*}
For such interfaces, the following energy identity can be verified integrating by parts:
\begin{align}\label{eq:EnergyId}
&\calE(\Sigma_2;w,\btau)=
\calE(\Sigma_1;w,\btau)-\int_{\Gamma_{\Sigma_1,\Sigma_2}} w\,\btau\cdot\bn^x_\Omega\di S\\
&+\int_{Q_{\Sigma_1,\Sigma_2}}\bigg(
\Big(\nabla w+\der{\btau}t \Big) \cdot \btau
+\Big(\nabla\cdot\btau+\frac{1}{c^2}\der{w}t\Big) w\bigg)
\di V.
\nonumber
\end{align}
If $(v,\bsigma)$ is the solution of IBVP~\eqref{eq:IBVP}, then we have
$$\calE(\Sigma_2;v,\bsigma)=
\calE(\Sigma_1;v,\bsigma)-\int_{\Gamma_{\Sigma_1,\Sigma_2}} v\,\bsigma\cdot\bn^x_\Omega\di S.
$$
If $g_D=g_N=g_R=0$ in their domains of definition, since $\vartheta\ge0$, we have 
$\calE(\Sigma_2;v,\bsigma)\le\calE(\Sigma_1;v,\bsigma)$, i.e.\ energy is dissipated.
If moreover $\FR=\emptyset$, then $\calE(\Sigma_2;v,\bsigma)=\calE(\Sigma_1;v,\bsigma)$, i.e.\ energy is preserved.

Expanding $\ell(\vhp,\Shp)\ge\Tnorm{(\vhp,\Shp)}\DG^2$ as in \cite[eq.~(28)]{SpaceTimeTDG}, it is easy to prove that the method \eqref{eq:TDG} is \emph{dissipative} and that energy is dissipated by the discrete solution jumps across mesh interfaces: if $g_D=g_N=g_R=0$, the energy of the discrete solution $(\vhp,\Shp)$ of the Trefftz-DG formulation satisfies
$$\mmbox{\calE(T;\vhp,\Shp)\le 
\calE(0;v_0,\bsigma_0)
.}$$ 
Using $|\bn^x_F|\le\gamma c^{-1}n^t_F$ on $\Fspa$ and the weighted Cauchy--Schwarz inequality, we have lower and upper bounds for the energy:
if $\Sigma\subset\Fspa\cup\FO\cup\FT$ is a space-like interface composed by element faces, then for all $ (w,\btau)\in L^2(\Sigma)^{1+n}$
\begin{align}\label{eq:EnergyLowerBound}
&\frac12\int_\Sigma (1-\gamma) n^t_\Sigma\Big(\frac{w^2}{c^2}+|\btau|^2\Big)\di S
\quad\le\;
\calE(\Sigma; w, \btau)
\\&
\le\frac12\int_\Sigma (1+\gamma) n^t_\Sigma\Big(\frac{w^2}{c^2}+|\btau|^2\Big)\di S
.\nonumber
\end{align}


\subsection{Stability of the auxiliary problem}\label{s:Stability}

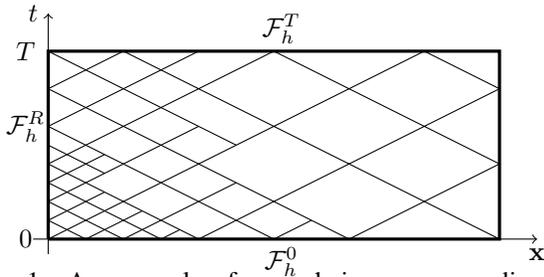
\begin{figure}[H]
\centering
\begin{tikzpicture}[scale=1]
\draw [->] (0,-0.2)--(0,3);\draw(-0.2,3)node{$t$};
\draw [->] (-0.2,0) -- (6.5,0);\draw(6.5,-0.2)node{$\bx$}; 
\draw [very thick] (0,0)--(6,0)--(6,2.5)--(0,2.5)--(0,0);
\draw (0,0)--(5,2.5)--(6,2)--(2,0)--(0,1)--(3,2.5)--(6,1)--(4,0)--(0,2)--(1,2.5)--(6,0);
\draw (3.5,0.25)--(3,0)--(0,1.5)--(2,2.5);
\draw (2.5,0.75)--(1,0)--(0,0.5)--(2,1.5); 
\draw (0,2.5)--(2.5,1.25);
\draw (1.75,0.125)--(1.5,0)--(0,0.75)--(0.75,1.125);
\draw (1.25,0.375)--(0.5,0)--(0,0.25)--(0.75,0.625);
\draw (0.75,0.875)--(0,1.25);
\draw(3.1,-0.3)node{$\FO$}; 
\draw(3.1,2.8)node{$\FT$}; 
\draw(-0.3,1.5)node{$\FR$}; 
\draw(-0.3,0)node{$0$}; 
\draw(-0.3,2.5)node{$T$}; 
\end{tikzpicture}\vspace{-7mm}
\caption{An example of a mesh in one space dimension ($n=1$) satisfying the assumptions of Proposition \ref{prop:DualStability}.
All internal mesh faces are space-like.
Not all mesh elements are 2-simplices (triangles), but the faces are 1-simplices (segments), so all integrals in \eqref{eq:TDG} are easy to compute.
For this mesh, the parameter $N$ in the proof of Proposition \ref{prop:DualStability} is equal to 16.
\label{fig:Mesh}
}
\end{figure}

We are now ready to prove bound \eqref{eq:DualStability} which guarantees the error bound in $L^2(Q)$ norm \eqref{eq:L2Error}.
We note that, under assumption (ii) below, the wave speed $c$ must be constant throughout $Q$.
A mesh satisfying this assumption is depicted in Figure \ref{fig:Mesh}.

\begin{proposition}\label{prop:DualStability}
Consider the IBVP \eqref{eq:zzIBVP} and assume that
\begin{itemize}
\item[(i)]
$\GD=\GN=\emptyset$, i.e.\ only Robin boundary conditions are allowed, $\deO=\GR$; and
\item[(ii)]
$\Ftime=\emptyset$, i.e.\ no time-like mesh interfaces are allowed.
\end{itemize}
Then the stability estimate \eqref{eq:DualStability} holds true with
\begin{align}
\Cstab^2\!=\!2T\bigg(
\!N\Norm{\frac{4(1+\gamma^2)}{(1-\gamma)^2}}_{L^\infty(\Fspa)}
\!\!+\Norm{\frac1{\delta(1-\delta)}}_{L^\infty(\FR)}
\!\bigg)
\nonumber
\end{align}
where $N$ is the minimal number of space-like interfaces $\Sigma_1,\ldots,\Sigma_N$ such that 
$\Fspa\subset\bigcup_{1\le j\le N-1}\Sigma_j$ and $0\le f_{\Sigma_1}\le\cdots\le f_{\Sigma_{N-1}}\le f_{\Sigma_N}=T$.

\end{proposition}
\textit{Proof.}\;
Applying the energy identity \eqref{eq:EnergyId} to the solution $(z,\bzeta)$ of the IBVP \eqref{eq:zzIBVP}, we have that for any two space-like interfaces $\Sigma_1,\Sigma_2$ with $f_{\Sigma_1}\le f_{\Sigma_2}$,
\begin{align}\label{eq:zzEnergyEvolution}
\calE(\Sigma_2; z,\bzeta)\le
\calE(\Sigma_1; z,\bzeta)
+\int_{Q_{\Sigma_1,\Sigma_2}}\big(\bPhi \cdot \bzeta+\psi z\big)\di V
\end{align}
(equality holds if $\Gamma_{\Sigma_1,\Sigma_2}\cap\FR$ has vanishing $n$-dimensional measure).
This implies a bound in space--time $L^2$ norm:
\begin{align*}
&\Norm{c^{-1}z}^2_{L^2(Q)}+\Norm{\bzeta}^2_{L^2(Q)^n}
=2\int_0^T \calE(t;z,\bzeta)\di t
\\&
\overset{\eqref{eq:zzEnergyEvolution}}\le 2\int_0^T 
\Big(\underbrace{\calE(0;z,\bzeta)}_{=0} 
+\int_\Omega\int_0^t(\bPhi\cdot\bzeta+\psi z)\di s\di \bx\Big)\di t\\
&\le 2T \Big(\Norm{c^{-1}z}^2_{L^2(Q)}+\Norm{\bzeta}^2_{L^2(Q)^n}\Big)^{1/2}
\\&
\qquad\cdot\Big(\Norm{c\psi}^2_{L^2(Q)}+\Norm{\bPhi}^2_{L^2(Q)^n}\Big)^{1/2},
\end{align*}
from which
\begin{align}\label{eq:zzL2QBound}
\Norm{c^{-1}z}^2_{L^2(Q)}\!+\Norm{\bzeta}^2_{L^2(Q)^n}
\!\le 4T^2 \big(\Norm{c\psi}^2_{L^2(Q)}\!+\!\Norm{\bPhi}^2_{L^2(Q)^n}\big).
\end{align}
For every space-like mesh interface $\Sigma\subset\Fh$ we control the corresponding term in \eqref{eq:DualStability} with the energy term:
\begin{align}
\label{eq:zzCgamma}
&2\Norm{\Big(\frac{(1\!+\!\gamma^2)n^t_\Sigma}{1-\gamma}\Big)^{1/2}\frac zc}^2_{L^2(\Sigma)}
\!+2\Norm{\Big(\frac{(1\!+\!\gamma^2)n^t_\Sigma}{1-\gamma}\Big)^{1/2}\bzeta}^2_{L^2(\Sigma)^n}
\\
&\le \underbrace{\Norm{\frac{4(1+\gamma^2)}{(1-\gamma)^2}}_{L^\infty(\Fspa)}}_{=:C_\gamma}
\frac12\int_\Sigma (1-\gamma) n^t_\Sigma\Big(\frac{z^2}{c^2}+|\bzeta|^2\Big)\di S
\nonumber\\&
\overset{\eqref{eq:EnergyLowerBound}}\le C_\gamma\;\calE(\Sigma;z,\bzeta).
\nonumber
\end{align}
We partition the faces in $\Fspa$ into $(N-1)$ interfaces $\Fspa=\bigcup_{1\le j\le N-1}\Sigma_j$  such that 
$0\le f_{\Sigma_1}\le\cdots\le f_{\Sigma_{N-1}}\le 
T$ and denote $\Sigma_N=\Omega\times\{T\}$.
We now control all terms on the space-like faces:
\begin{align*}
&2\Norm{\Big(\frac{(1+\gamma^2)n^t_F}{1-\gamma}\Big)^{1/2}\frac zc}^2_{L^2(\Fspa\cup\FT)}\\
&\qquad+2\Norm{\Big(\frac{(1+\gamma^2)n^t_F}{1-\gamma}\Big)^{1/2}\bzeta}^2_{L^2(\Fspa\cup\FT)^n}
\\&
\overset{\eqref{eq:zzCgamma}}\le C_\gamma \sum_{j=1}^N \calE(\Sigma_j;z,\bzeta)
\\&
\overset{\eqref{eq:zzEnergyEvolution}}\le C_\gamma \sum_{j=1}^N \int_{Q_{\Omega\times\{0\},\Sigma_j}}\big(\bPhi \cdot \bzeta+\psi z\big)\di V
\\&
\overset{\eqref{eq:zzL2QBound}}\le
2C_\gamma N T\big(\Norm{c\psi}^2_{L^2(Q)}+\Norm{\bPhi}^2_{L^2(Q)^n}\big).
\end{align*}
From assumptions (i) and (ii), in \eqref{eq:DualStability} there are no terms on $\Ftime,\FD$ and $\FN$, so we are now left with the terms on $\FR$:
using the Robin boundary condition $\frac\tht c z=\bzeta\cdot\bn_\Omega^x$,
the energy identity \eqref{eq:EnergyId} and the $L^2(Q)$ stability bound \eqref{eq:zzL2QBound}, we have
\begin{align*}
&\Norm{\Big(\frac c{(1-\delta)\tht}\Big)^{1/2}\bzeta\cdot\bn^x_\Omega}^2_{L^2(\FR)}
+\Norm{\Big(\frac \tht{\delta c}\Big)^{1/2}z}^2_{L^2(\FR)}
\\&
\le\underbrace{
\Norm{\frac1{\delta(1-\delta)}}_{L^\infty(\FR)}}_{=:C_\delta}
\int_\FR\!\Big(
\delta\frac c\tht(\bzeta\cdot\bn^x_\Omega)^2+ (1-\delta)\frac\tht cz^2\Big)\di S
\\&
\overset{\tht z=c\bzeta\cdot\bn_\Omega^x}=
C_\delta
\int_\FR z\:\bzeta\cdot\bn^x_\Omega\di S
\\&
\overset{\eqref{eq:EnergyId}}= C_\delta\bigg(
\underbrace{\calE(0;z,\bzeta)}_{=0}-\underbrace{\calE(T;z,\bzeta)}_{\ge0}
+\int_Q\big(\bPhi\cdot\bzeta+\psi\, z\big)\di V
\bigg)
\\&
\overset{\eqref{eq:zzL2QBound}}\le C_\delta 2T \big(\Norm{c\psi}^2_{L^2(Q)}+\Norm{\bPhi}^2_{L^2(Q)^n}\big).
\end{align*}
Combining this inequality with the previous one we obtain the assertion with $\Cstab^2=2T(C_\gamma N+C_\delta)$.
\hfill$\square$

Assumption (ii) in Proposition \ref{prop:DualStability} requires that all the internal mesh faces are space-like; Figure \ref{fig:Mesh} shows a mesh of this kind.
The meshes that satisfy this condition allow the Trefftz-DG method to be treated as a ``semi-explicit'' scheme as in \cite{FalkRichter1999,MoRi05}:
if the elements are suitably ordered, the discrete solution can be computed sequentially solving a local problem for each element.
This also allows a high degree of parallelism.
If the ``\emph{tent-pitching}'' algorithm of \cite{EGSU05} 
 is used to construct the mesh and the ``macro elements'' of \cite{MoRi05} are taken as elements, then the mesh obtained satisfies the assumptions of Proposition \ref{prop:DualStability}. 
The fact that the elements obtained in this way do not have simple shapes such as $(n+1)$-simplices is not a computational difficulty: all integrals in the formulation \eqref{eq:TDG} are defined on mesh faces, which are $n$-simplices, thus no quadrature on complicated shapes needs to be performed.
This is due to the Trefftz property, so this advantage is not available to discretisations employing standard (non-Trefftz) local spaces.

\section{
Polynomial Trefftz spaces}\label{s:Spaces}

So far we have not specified any discrete space $\bV \Th$: the only condition we imposed is the Trefftz property, i.e.\ $\bV \Th\subset\bT\Th$.
Given $p\in\N$, a simple choice is to define  $\bV \Th:=\prod_{K\in\calT_h}\bV_p(K)$, where $\bV_p(K):=\bT(K)\cap\mathbb{P}^p(\R^{n+1})^{1+n}$ is the space of the solutions $(w,\btau)$ of the wave equation in $K$ that are polynomials of degree at most $p$.

For high polynomial degree $p$, the Trefftz space $\bV_p(K)$ has much lower dimension than the space-time, vector-valued, full polynomial space $\mathbb{P}^p(\R^{n+1})^{1+n}$:
\begin{align*}
&\dim\bV_p(K)=\binom{p+n}n\frac{2p+n+2}{p+1}-1=\mathcal{O}_{p\to\infty}(p^{n}),\\
&\dim\mathbb{P}^p(\R^{n+1})^{1+n}
=\binom{p+n}n(p+n+1) = \mathcal{O}_{p\to\infty}(p^{n+1}).
\end{align*}
In dimension $n=1$, it has been proved in \cite[\S5.3]{SpaceTimeTDG} that these two spaces have the same orders of approximation in the mesh-width $h_K:=\mathrm{diam}\, K$ and in the polynomial degree $p$.
This is true both for solutions with bounded Sobolev regularity (algebraic orders in $h_K$ and $p$) and for analytic solutions (exponential orders in $p$).
Thus the order of convergence in terms of numbers of degrees of freedom can be considerably higher for   the Trefftz-DG method compared to similar DG schemes based on standard spaces (see e.g.\ \cite[Fig.~2]{SpaceTimeTDG}).
Similar results have been observed in higher dimensions \cite{EKSW2,WTF14}, but the proof of the orders of convergence is open. 

A basis of $\bV_p(K)$ can be con\-struc\-ted using ``transport polynomials'', or ``polynomial waves'', in the form $(\der{}tP_{\ell,j},-\nabla P_{\ell,j})$ with $P_{\ell,j}(\bx,t)=(\bx\cdot\bd_j-ct)^\ell$ for $0\le\ell\le p$, where the propagation directions $\bd_j\in\R^n$, $|\bd_j|=1$,  are suitably chosen.
Non-polynomial Trefftz basis functions can easily be constructed, for example, from $f(\bx\cdot\bd_j-ct)$ for any smooth function $f:\R\to\R$.

\label{ultima-pagina:}\end{multicols}
\end{document}